\documentclass[authoryear,preprint,review,10pt]{elsarticle}
\usepackage{subfigure}
\usepackage[figuresright]{rotating}
\usepackage{bm}
\usepackage{color}
\RequirePackage{amsthm,amsmath,amsfonts,amssymb}
\theoremstyle{plain}

\numberwithin{equation}{section}

\journal{}

\begin{document}
	
	\begin{frontmatter}
		
		\title{Modile as a conservative tail risk measurer: the solution of an optimisation problem with  0-1 loss function }
		
		\author[mymainaddress]{Keming Yu\corref{mycorrespondingauthor}}
		\cortext[mycorrespondingauthor]{Corresponding author}
		\ead{keming.yu@brunel.ac.uk}
		
		\author[mysecondaryaddress]{Rong Jiang}
		\ead{jiangrong@sspu.edu.cn}
		
		\author[mythirdaryaddress]{Chi Tim Ng}
		\ead{timothyng@hsu.edu.hk}

		\address[mymainaddress]{Brunel University London, UB8 3PH, UK}
		\address[mysecondaryaddress]{Shanghai Polytechnic University, Shanghai, 201620, China}
		\address[mythirdaryaddress]{The Hang Seng University of Hong Kong, 	Hang Shin Link, Siu Lek Yuen, Sha Tin, Hong Kong, China}

		\begin{abstract}
			Quantiles and expectiles, which are two important concepts and tools in tail risk measurements, can be regarded as an extension of median and mean, respectively.  Both of these tail risk measurers can actually be embedded in a common framework of $L_p$ optimization with the absolute loss function ($p=1$) and quadratic loss function ($p=2$), respectively. When 0-1 loss function is frequently used in statistics, machine learning and decision theory, this paper introduces an 0-1 loss function based $L_0$ optimisation problem for  tail risk measure and names its solution as {\it modile}, which can be regarded as an extension of mode.
Mode, as another measure of central tendency, is more robust than expectiles with outliers and  easy to compute than quantiles.
However, mode based extension for tail risk measure
			is new. This paper shows that the proposed {\it modiles} are not only  more  conservative than quantiles and expectiles
for skewed and heavy-tailed distributions, but also providing or including  the unique interpretation of  these  measures.
Further, the {\it modiles}  can be regarded as a type of   generalized quantiles and doubly truncated tail measure whcih have recently attracted a lot of attention in the literature.
The asymptotic properties of the corresponding sample-based estimators of {\it modiles} are provided, which, together with numerical analysis results,
show that the proposed {\it modiles} are promising for tail  measurement.
			
		\end{abstract}
		
		
	\end{frontmatter}

	\section{Introduction}
	Quantile \citep{r1} and expectile \citep{r2} have been attracted considerable attention in the literature with wide application in statistics, economics, finance and risk management,
	see \cite{breckling1988m}, \cite{waltrup2015expectile}, \cite{r4}, \cite{r5}, \cite{r3},  \cite{madrid2022risk}, \cite{dimitriadis2022realized}, \cite{dimitriadis2022characterizing} and among others.
	\par
	For $\tau\in (0,1)$,
	the $\tau$th quantile $q_{\tau}$ of a random variable $Y$ can be regarded as a skew modification of median as it minimises an asymmetric absolute loss function:
	\begin{equation}
		\begin{split}
			q_{\tau}=\arg\min_{\theta\in \mathbb{R}}
			\text{E}\left\{ \left|\tau-I(Y\leq \theta)  \right|\cdot|Y-\theta|\right\}.
		\end{split}
	\end{equation}
	This is, $q_{\tau}$ is the solution of a $L_1$ optimisation problem. When $\tau=1/2$, the quantile of $Y$ reduces to its mean.
	\par
	The $\tau$th expectile $\xi_{\tau}$ of a random variable $Y$ can be regarded as a skew modification of mean as it minimises an asymmetric quadratic loss function:
	\begin{equation}
		\xi_{\tau}=\arg\min_{\theta\in \mathbb{R}} \text{E}\left\{|\tau-I(Y\leq \theta)|\cdot(Y- \theta)^2\right\}.
	\end{equation}
	This is,  $\xi_{\tau}$ is the solution of  a  $L_2$ optimisation problem. When $\tau=1/2$, the expectile of $Y$ reduces to its mean.
	Note that $q_{\tau}$ equivalently (\citep{r2})
	satisfies
	\begin{equation}
		\frac{P(Y\geq q_{\tau})}{P(Y<q_{\tau})}=
		\frac{1-F_Y(q_{\tau})}{F_Y(q_{\tau})}
		=\frac{1-\tau}{\tau},
	\end{equation}
	and $\xi_{\tau}$ equivalently
	satisfies
	\begin{equation}
		\frac{\int_{\xi_{\tau}}^{+\infty}(y-\xi_{\tau})\,dF_Y(y)}{\int_{-\infty}^{\xi_{\tau}}(\xi_{\tau}-y)\,dF_Y(y)}
		=\frac{E(Y-\xi_{\tau})_+}{E(Y-\xi_{\tau})_-}=\frac{1-\tau}{\tau},
	\end{equation}
	where $F_Y(\cdot)$ is the cumulative distribution function (cdf) of a random variable $Y$.
	Here $(1-\tau)/{\tau}$ is the targeted probability ratio and the expectile level ratio or targeted gain-loss ratio for (1.3) and (1.4), respectively.
	
	\par
	Mode \citep{r8}, as another  most important characteristic of a random variable, usually represents the value where probability
	density is the biggest for a continuous random variable, and the most probable value for a discrete random variable but links to tail  risk measure rarely.
	However, based on an optimisation problem with an 0-1 loss function, this paper introduces a new tail risk measurement, which can  be an extension of mode and
	include the unique interpretation of both quantiles and expectiles.
	\par
	The remainder of this paper is organized as follows. In Section 2, the definition of modile is proposed.
	Section 3 shows that  modiles are more conservative than the popular quantiles and expectiles for skewed and  heavy-tailed distributions.
The estimation algorithm and asymptotic property of modiles are developed in Section 4.
	Both simulation examples and the application of modiles on real data
	are given in Section 5.  Some technical proofs are provided in the Appendix.
	
	\section{The Modile and its basic properties}
	
	Note that the mode $\nu$ of a random variable $Y$ could be defined in terms of the $L_0$ optimisation problem
 $\nu=\arg\min_{\theta\in \mathbb{R}} \text{E}\left\{L(Y,\theta) \right\}$ with the 0-1 loss unction $L(Y,\theta)=I(|Y-\theta|>\delta)$.
Where  $\delta>0$ is a constant. See \citep{r6, r7} and among others.
	\par
	Also note that the defined mode above often requires the density function $f_Y(y)$ of $Y$ is strictly unimodal.
	To relax the restriction and include  non-monotone density, asymmetric density and even multi-modual, we may consider
	an alternative 0-1 loss function with a lower limit $-h_1$ and  upper limit $h_2$  and then formulate  the $L_0$ optimisation problem as
	\[\nu=\arg\min_{\nu\in \mathbb{R}} \text{E}\{I(-h_1\leq Y-\theta \leq h_2) \}.\]
	In fact, this 0-1 loss function should define  $\nu$  as the mode of  the random variable $Y$ \citep{r11}. When $h_1=h_2$ and $\tau=0.5$, modile (2.1) is equal to mode \citep{r6}.
	
	Then, along the same line of defining quantiles and expectiles as the  optimal predictors
	under an asymmetric absolute loss function and asymmetric least squares loss function respectively, we define the {\it  modile} as an asymmetric
	0-1 loss function based $L_0$ optimisation as:
	\begin{equation}
		\begin{split}
						\nu_{\tau}=\arg\max_{\theta\in \mathbb{R}} \text{E}\left\{|\tau-I(Y\leq \theta)|\cdot I(-h_1\leq Y-\theta \leq h_2)\right\}
\end{split}
	\end{equation}
or
\begin{equation}
\begin{split}
			\nu_{\tau}=\arg\min_{\theta\in \mathbb{R}} \text{E}\left[|\tau-I(Y\leq \theta)|\cdot \left\{I(Y-\theta<-h_1)+I(Y-\theta>h_2)\right\}\right],
		\end{split}
	\end{equation}
	for two fixed positive number $h_1$ and $h_2$. That is,  $\nu_{\tau}$ is defined as the $\tau$th modile of a random variable $Y$ for any $0<\tau<1$.
	Based on the equations (2.1) or (2.2), we have
	\begin{equation*}
\begin{split}
					\text{E}\left\{\rho_{\tau}(Y-\nu_{\tau},h_1,h_2) \right\}
			=&\tau\,P(Y>\nu_{\tau}+h_2)+(1-\tau)\,P(Y<\nu_{\tau}-h_1)\\
			=&\tau-\tau\, F_Y(\nu_{\tau}+h_2)+(1-\tau)\,F_Y(\nu_{\tau}-h_1),
		\end{split}
	\end{equation*}
	where
	$$\rho_{\tau}(u,h_1,h_2)=|\tau-I(u<0)|\cdot \left\{I(u<-h_1)+I(u>h_2)\right\}=
	\begin{cases}
		\tau, & \text{if } u > h_2\\
		1-\tau, & \text{if } u<-h_1\\
		0, & \text{otherwise. }
	\end{cases}$$
	So, we have
	\begin{equation}
		\begin{split}
			\frac{d\,E\rho_{\tau}(Y-\nu_{\tau},h_1,h_2)}{d\,\nu_{\tau} }=-\tau\, f_Y(\nu_{\tau}+h_2)+(1-\tau)\,f_Y(\nu_{\tau}-h_1)=0,
		\end{split}
	\end{equation}
	where $f_Y(\cdot)$ is the probability density function (pdf) of a random variable $Y$.
	We obtain
	\begin{equation}
		\frac{f_Y(\nu_{\tau}+h_2)}{f_Y(\nu_{\tau}-h_1)}=\frac{1-\tau}{\tau}.
	\end{equation}
	Here $f_Y(\nu_{\tau}+h_2)$ and $f_Y(\nu_{\tau}-h_1)$ can be interpreted as the most likely gain and most likely loss of $Y$, respectively.
	\par
	In the following, in some special cases of $f_Y(\cdot)$, let's specify the modile $\nu_{\tau}$ via (2.4).
	\par
	(1) If $Y\sim Normal(\mu,\sigma^2)$, we have $\nu_{\tau}=\mu+\frac{h_1-h_2}{2}-\frac{\sigma^2}{h_1+h_2}\ln\left(\frac{1-\tau}{\tau}\right)$.
Let  $h_1=h_2$, then $\nu_{\tau}=\mu$ ($\mu$ is the mode) under $\tau=0.5$.
	\par
	(2) If $Y\sim Laplace(\mu,\lambda)$, $\nu_{\tau}=\mu+(h_1-h_2)/2-\frac{\lambda}{2}\ln\left(\frac{1-\tau}{\tau}\right)$ for $\min\{h_1,h_2\}\geq |\nu_{\tau}-\mu|$.
	If $h_1=h_2$, $\nu_{\tau}=\mu$ ($\mu$ is the mode) under $\tau=0.5$.
	\par
	(3) If $Y\sim Gamma(\alpha,\beta)$, its density function is $f(x)=\frac{\beta^{\alpha}}{\Gamma(\alpha)}x^{\alpha-1}e^{-\beta x}$ with $x\geq0,\alpha>0,\beta>0$.
	We have $\nu_{\tau}=(Ah_1+h_2)/(A-1)$,
	where $A=\frac{1-\tau}{\tau}\exp\{\frac{\beta}{\alpha-1}(h_1+h_2)\}$. When $h_1$ and $h_2$ satisfy
	$
	h_1+h_2=B\ln\left(\frac{B+h_2}{B-h_1}\right),
	$
	where $B=(\alpha-1)/\beta$ is the mode, then $\nu_{\tau}=B$ under $\tau=0.5$.
	\par
	(4) If $f_Y(x)=12(b-a)^3\cdot(x-\frac{a+b}{2})^2\cdot I(a\leq x \leq b)$, the mode is $a$ and $b$, we have
	$
	\nu_{\tau}=\left\{{\frac{a+b}{2}-h_2+(\frac{a+b}{2}+h_1)A}\right\}/(1+A),
	$
	where $A=\pm\sqrt{\frac{1-\tau}{\tau}}$. When $\tau=0.5$, $\nu_{\tau}=(a+b-h_2+h_1)/2=a$ (or $b$) under $h_1=a$, $h_2=b$ (or $h_1=b$, $h_2=a$).
	\par

\textcolor{black}{In the above examples, the modile $\nu_{\tau}$ has a close-form expression under some density functions. Moreover, mode is a special case of modile for symmetrically distributed loss particularly when $h_1=h_2$. In the following example, we set $h_1=a+|b-c|$ and $h_2=a+|b+c|$, where $a,b,c$ are standard deviation, mean and skewness of data, respectively. The numberical results of the quantiles, expectiles and modiles for Normal distribution, Laplace distribution and Gamma distribution are shown in Table 1 and Figures 1-3.}

	\begin{table}
		\footnotesize
		\caption{The modile, quantile and expectile
			under different errors and $\tau$s. }
		\centering
		\begin{tabular}{@{}r|ccc|ccc|ccc@{}}
			\hline
			$\bm{Y}$&&Normal(0,1)&&&Laplace(1,2)&&&Gamma(8,7)\\
			\hline
			$\tau$&modile&quantile&expectile&modile&quantile&expectile&modile&quantile&expectile\\
			\hline
			0.1&-1.099&-1.282&-0.861&-1.197&-2.219&-1.404&0.855&0.665&0.833\\
			0.2&-0.693&-0.842&-0.549&-0.386&-0.833&-0.452&0.875&0.797&0.938\\
			0.3&-0.424&-0.524&-0.337& 0.153&-0.022& 0.135&0.901&0.902&1.014\\
			0.4&-0.203&-0.253&-0.162& 0.595& 0.554& 0.592&0.936&0.999&1.080\\
			0.5& 0.000& 0.000& 0.000& 1.000& 1.000& 1.000&0.987&1.096&1.143\\
			0.6& 0.203& 0.253& 0.162& 1.405& 1.446& 1.408&1.066&1.199&1.209\\
			0.7& 0.424& 0.524& 0.337& 1.847& 2.022& 1.865&1.206&1.316&1.283\\
			0.8& 0.693& 0.842& 0.549& 2.386& 2.833& 2.452&1.525&1.462&1.377\\
			0.9& 1.099& 1.282& 0.861& 3.197& 4.219& 3.404&2.972&1.682&1.523\\
			\hline
		\end{tabular}
	\end{table}
	\begin{figure}
		\centering
		{\label{figa} 
			\includegraphics[width=4.5in]{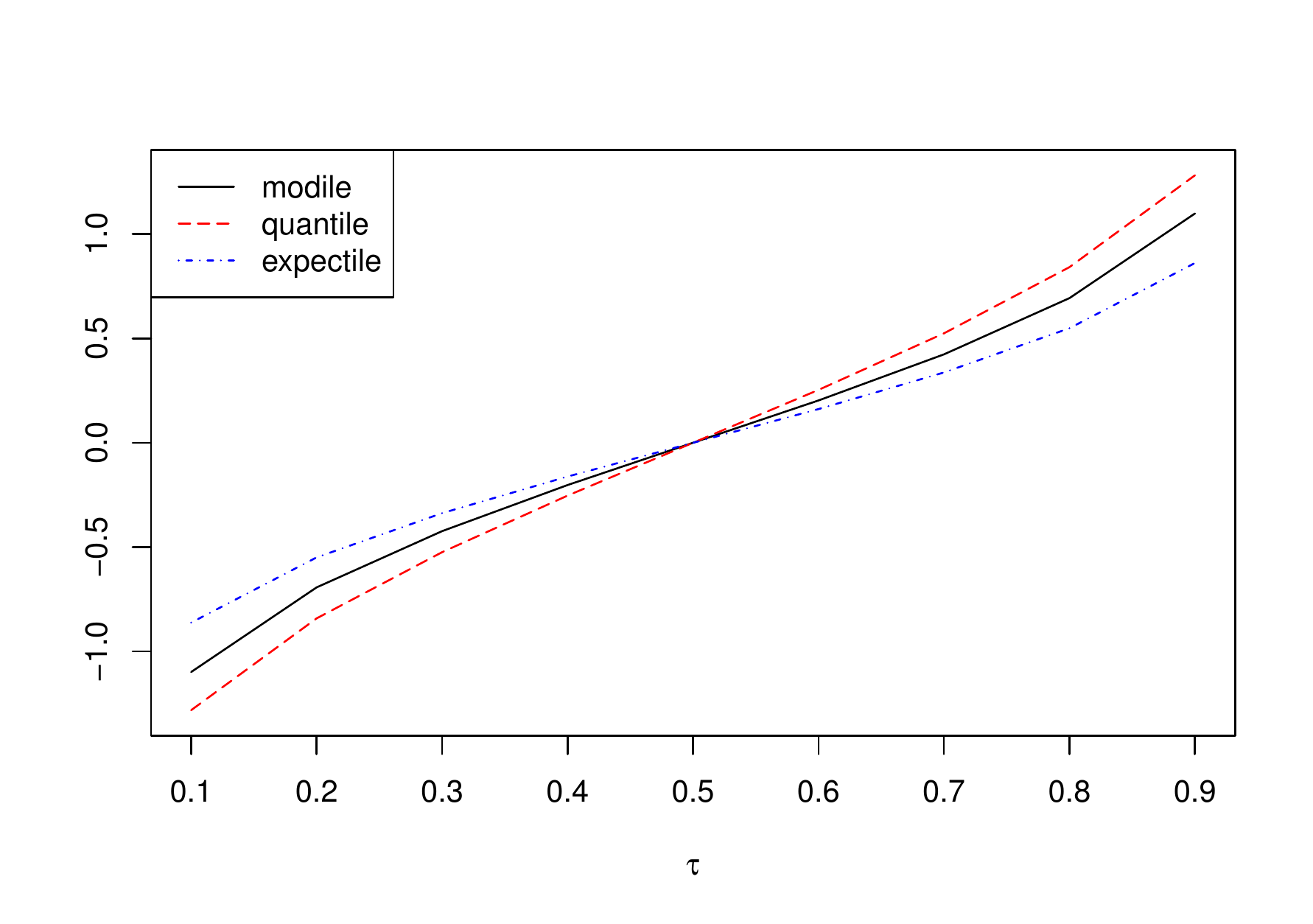}}
		\caption{The modile, quantile and expectile for $\bm{Y}\sim$Normal(0,1).}
	\end{figure}
	
	\begin{figure}
		\centering
		{\label{figa} 
			\includegraphics[width=4.5in]{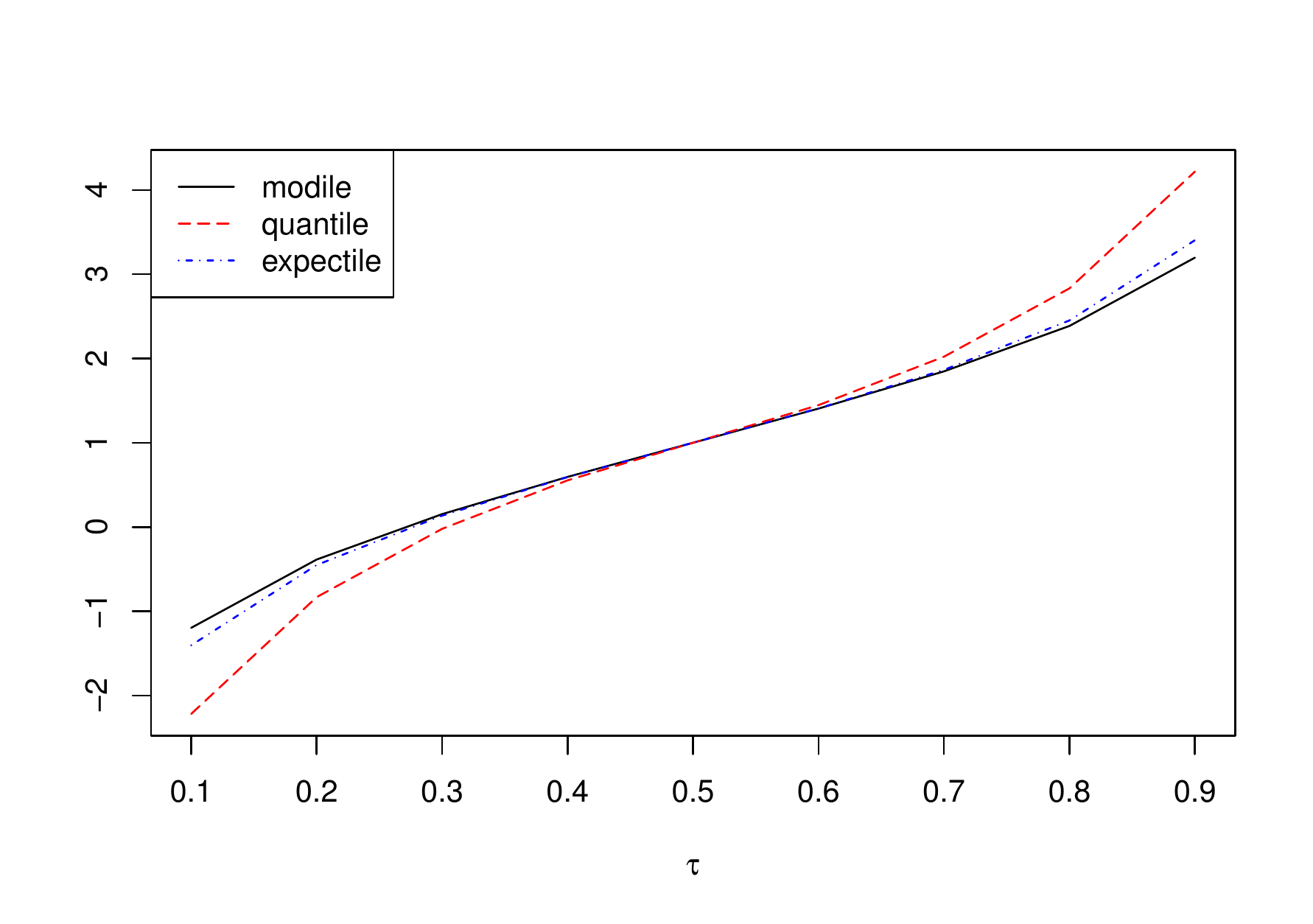}}
		\caption{The modile, quantile and expectile for $\bm{Y}\sim$Laplace(1,2).}
	\end{figure}
	
	\begin{figure}
		\centering
		{\label{figa} 
			\includegraphics[width=4.5in]{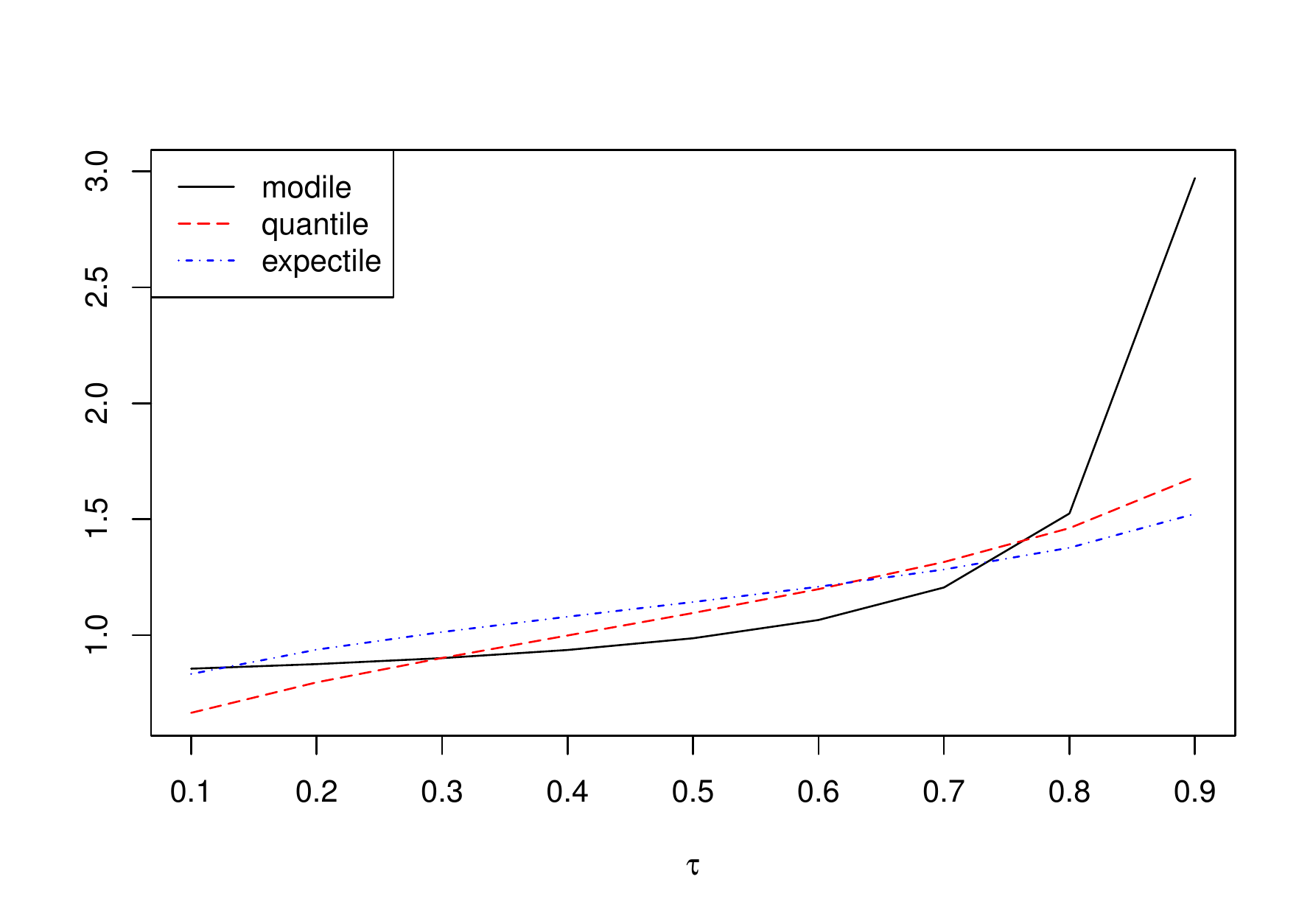}}
		\caption{The modile, quantile and expectile for $\bm{Y}\sim$Gamma(8,7).}
	\end{figure}
	
\newpage
	{\bf Theorem 1.} Let $\delta=h_1=h_2 \equiv \delta$, for any $x\in(0,\delta)$, the $\nu_{\tau}$  satisfies that
	\begin{equation}
		\frac{f_Y(\nu_{\tau}+x)}{f_Y(\nu_{\tau}-x)}=\frac{1-\tau}{\tau}.
	\end{equation}
	Then for any $\delta$, the following equation holds:
	\begin{equation}
		\begin{split}
			\frac{P(\nu_{\tau} \leq Y \leq \nu_{\tau}+\delta)}{P(\nu_{\tau}-\delta \leq Y \leq \nu_{\tau})}=\frac{1-\tau}{\tau},
		\end{split}
	\end{equation}
	and
	\begin{equation}
		\frac{E\{(Y-\nu_{\tau}) \,I(\nu_{\tau} \leq Y \leq \nu_{\tau}+\delta)\}}{E\{(\nu_{\tau}-Y)\,I(\nu_{\tau}-\delta \leq Y \leq \nu_{\tau})\}}=\frac{1-\tau}{\tau}.
	\end{equation}

	When $\delta \rightarrow \infty$, we have
	\begin{equation}
		\frac{P(Y \geq \nu_{\tau})}{P(Y \leq \nu_{\tau})}=\frac{1-\tau}{\tau},
	\end{equation}
	and
	\begin{equation}
		\frac{E\{(Y-\nu_{\tau}) \,I( Y \geq \nu_{\tau})\}}{E\{(\nu_{\tau}-Y)\,I(Y \leq \nu_{\tau})\}}=\frac{1-\tau}{\tau}.
	\end{equation}
	\\
	\par
	Theorem 1 shows that modile of $Y$ are determined by its tail range (interval)  probabilities.
	The left side of the equation (2.6) can be interpreted as the  probability ratio of the ranges of gain and loss, and the right side $(1-\tau)/{\tau}$
is the targeted most-likely gain-loss range ratio.
	(2.8) shows that quantiles are the special cases of modile.
	Furthermore, we have the expectations of the tail range of $Y$ in (2.6),
	and the expectiles of $Y$ are also the special cases of modile as when $\delta \rightarrow \infty$ in (2.8).
	\par
	{\bf Remark 1:}
	 For any $\tau \in (0,1)$, the following density function $f_Y(\cdot)$ satisfies the condition (2.5),
	\begin{equation*}
		f_Y(x)=
		\begin{cases}
			~2(1-\tau)f(x), & \text{if } x\geq \nu\\\
			2\tau f(x), & \text{if } x<\nu,
		\end{cases}
	\end{equation*}
	where $f(\cdot)$ is a symmetric density function and $\nu$ is the mode of $f(\cdot)$. $f_Y(\cdot)$ is a function of $\tau$ because we make $\nu_{\tau}=\nu$.
	\par
	For example, $f(\cdot)$ is the density function of standard Normal distribution, then $\nu=\nu_{\tau}=0$. We have
	\begin{equation*}
		f_Y(x)=
		\begin{cases}
			\frac{2(1-\tau)}{\sqrt{2\pi}}\exp(-x^2/2), & \text{if } x\geq 0\\\
			\frac{2\tau}{\sqrt{2\pi}}\exp(-x^2/2), & \text{if } x<0,
		\end{cases}
	\end{equation*}
	satisfies the condition (2.4), and $\nu_{\tau}=q_{\tau}=\xi_{\tau}$.
	\par
	{\bf Remark 2:} Note that we may re-write
{\bf modiles} as the minimizer of
\begin{equation*}
		\nu_{\tau}=\arg\min_{\theta\in \mathbb{R}} \left[\tau \text{E}\{\Phi_1((Y-\theta)^+)\}+(1-\tau) \text{E}\{\Phi_2((Y-\theta)^-)\}\right],
	\end{equation*}
with convex loss function $\Phi_1(u)=I[u \leq Y\leq u+h_2]$ and $\Phi_2(u)=-I[u-h_1 \leq Y\leq u]$, which is  a type of  {\it Generalized quantiles} defined in
\citep{bellini2014generalized}  but without a need to satisfy the first-order condition of \citep{bellini2014generalized}.
The convex of $\Phi_i(\cdot)$ ($i=1, 2$) is due to the fact that the associated sets are convex.
\par
Equations (2.6)-(2.7) show that the proposes {\bf modiles} may also propose a type of doubly truncated tail probability,
 while   the doubly truncated tail conditional expectation $\text{E}(X|x_p< X < x_q)$  has recently attracted much attention in the literature
(\citep{roozegar2020moments,shushi2020multivariate}), where $x_p$ is the p-th quantile, and $0\leq p<q\leq 1.$

\section{Conservative risk measure for skewed and heavy-tailed distributions}

	For most studies in risk management, it has been found that
	Pareto distribution, which is a skewed and heavy-tailed distribution, describe quite well of the tail structure of actuarial and financial data.
Here we consider a Pareto distribution with shape parameter $\alpha\in(0,+\infty)$,
which is a continuous distribution on $(1,+\infty)$ with distribution function $F$ given by
	$F(y)=1-y^{-\alpha}.	$
	The probability density function $f(y)$ is given by
\begin{equation}
		f(y)=\alpha y^{-\alpha-1}.
	\end{equation}
	From the equation (3.1), the $\tau$-th modile $\nu_{\tau}$ is
	\begin{equation}
		\nu_{\tau}=\frac{h_1+Ah_2}{1-A},
		\end{equation}
	where $A=(1/\tau-1)^{\alpha+1}$. 	
	\par
It is easy to know that the quantile function of Pareto distribution is $q_{\tau}=(1-\tau)^{-1/\alpha}$. From the the result of Theorem 11 in \cite{bellini2014generalized},
the expectile function of Pareto distribution is $\xi_{\tau}\approx\{(1-\tau)(\alpha-1)\}^{-1/\alpha}$ under $\alpha>1$ and $\tau\to 1$.
Note that $\nu_{\tau}$ in (3.2), when $\tau\to 1$, we have $\nu_{\tau}=h_1<+\infty$, which is smaller than $q_{\tau}$ and $\xi_{\tau}$.
Therefore, modile is a more conservative risk measure than the quantile and expectile.
	\par
Generally, for the Pareto-like with tail index $\alpha$ distribution
\begin{equation*}
		F(y)=1-L(y) y^{-\alpha}
	\end{equation*}
for some function $L$ which is slowly varying at infinity.
We consider the $\tau$th quantiles $q_{\tau}$ conditional on the {\it modile}-range $[\nu_{\tau}-h_1, \,\, \nu_{\tau}+h_2]$:
 \begin{eqnarray*}
\tau&=&P[Y<q_{\tau}|\nu_{\tau}-h_1 \leq Y \leq \nu_{\tau}+h_2] \nonumber \\
&=&\frac{(\nu_{\tau}-h_1)^{-\alpha}-q_{\tau}^{-\alpha}}{(\nu_{\tau}-h_1)^{-\alpha}-(\nu_{\tau}+h_2)^{-\alpha}},
				\end{eqnarray*}
we have
\begin{equation*}
		q_{\tau}^{-\alpha}=\tau(\nu_{\tau}+h_2)^{-\alpha}+(1-\tau)\,(\nu_{\tau}-h_1)^{-\alpha}.
	\end{equation*}
Clearly, for any $\alpha\in(0,+\infty)$,  when $\tau \rightarrow 1$, $q_{\tau}=\nu_{\tau}+h_2>\nu_{\tau}$, and when $\tau \rightarrow 0$,
 $q_{\tau}=\nu_{\tau}-h_1<\nu_{\tau}$, so {\bf modiles} are   more generally conservative tail risk measures than {\bf quantiles}.
\par
Similar conclusion is true for {\bf expectiles}. In fact, consider the $\tau$th expectiles $\xi_{\tau}$ conditional on the {\it modile}-range $[\nu_{\tau}-h_1, \,\, \nu_{\tau}+h_2]$:
from
 \begin{equation*}
\frac{E(Y-\xi_{\tau})_+}{E(Y-\xi_{\tau})_-}=\frac{\int_{\nu_{\tau}}^{\xi_{\tau}+h_2}(x-\xi_{\tau})\,dF(x)}{\int_{\nu_{\tau}-h_1}^{\xi_{\tau}}(\xi_{\tau}-x)\,dF(x)   }=\frac{1-\tau}{\tau},
				\end{equation*}
and
 \begin{eqnarray*}
\int_{\nu_{\tau}-h_1}^{\xi_{\tau}}(\xi_{\tau}-x)\,dF(x)&=&\alpha \left\{ \frac{\xi_{\tau}^{1-\alpha}-(\nu_{\tau}-h_1)^{1-\alpha}}{1-\alpha}+ \xi_{\tau}\cdot\frac{\xi_{\tau}^{-\alpha}-(\nu_{\tau}-h_1)^{-\alpha}}{\alpha}\right\},  \nonumber \\
\int_{\nu_{\tau}}^{\xi_{\tau}+h_2}(x-\xi_{\tau})\,dF(x)&=&\alpha \left\{ \frac{(\xi_{\tau}+h_2)^{1-\alpha}-\xi_{\tau}^{1-\alpha}}{1-\alpha}+ \xi_{\tau}\cdot\frac{(\xi_{\tau}+h_2)^{-\alpha}-\xi_{\tau}^{-\alpha}}{\alpha}\right\}, \nonumber
				\end{eqnarray*}
under the Pareto-like heavy-tailed distribution, we have when $\tau \rightarrow 1$,
\begin{equation*}
\alpha (\nu_{\tau}+h_2)^{1-\alpha}+(1-\alpha) \xi (\nu_{\tau}+h_2)^{-\alpha}=\xi_{\tau}^{1-\alpha},
				\end{equation*}
which holds and only holds when $\xi_{\tau}=\nu_{\tau}+h_2$.
Along the same derivation, we can show that  when $\tau \rightarrow 0$, $\xi_{\tau}=\nu_{\tau}-h_1$.

	\section{Algorithm and asymptotic property of the estimation of modiles}
	Given a random sample $Y_1,\ldots,Y_n$ from $Y$, the estimator of $\hat{\nu}_{\tau}$ of $\hat{\nu}_{\tau}$
	 $\hat{\nu}_{\tau}$ can be obtained by the minimisation of the empirical  $L_0$  optimisation problem of the equation (2.2) under the 0-1 loss fcuntion as
\begin{equation}
	\hat{\nu}_{\tau}=\arg\min_{\theta\in \mathbb{R}}\left\{\tau-\frac{\tau}{n}\sum^n_{i=1}I(Y_i\leq \theta+h_2)+\frac{1-\tau}{n}\sum^n_{i=1}I(Y_i\leq \theta-h_1)\right\}.
\end{equation}
	We can use the multitask algorithm \citep{qin2020multitask} to find the numerical solutions of $\hat{\nu}_{\tau}$ easily.  See the R code attached.
	\par
	When the sample size $n$ tends to infty, the asymptotic distribution of $\hat{\nu}_{\tau}$ is given by Theorem 2.
	\\
	\par
	{\bf Theorem 2.} Assume that the density function of $Y$ is continuous differentiable, then we have
	$$
	n^{1/3}(\hat{\nu}_{\tau}-\nu_{\tau})\xrightarrow{d}\left[\frac{4\tau f_Y(\nu_{\tau}+h_2)}{\left\{\tau f'_Y(\nu_{\tau}+h_2)-(1-\tau)
		f'_Y(\nu_{\tau}-h_1)\right\}^{2}}\right]^{1/3}\hat{z},
	$$
	where $\xrightarrow{d}$ stands for convergence in the distribution, $\hat{z}=\arg\max_{z\in \mathbb{R}}\{Z(z)-z^2\}$, and $Z(z)$ is a two-sided Wiener-L\'{e}vy process through the origin with mean 0 and variance one per unit $z$.
	\par
	Note that the asymptotic distribution of $\hat{\nu}_{\tau}$ is the same as the estimation of the mode in \cite{r10} under $\tau=0.5$ and $h_1=h_2$.

	\section{Numerical studies}
	In this section, we first use Monte Carlo simulation studies to assess the finite sample performance of the proposed procedures and then demonstrate the application of the proposed methods with a real data analysis. All programs are written in \textsf{R} code.
	\subsection{Simulation example}
	In this section, we study the performance of the estimation of modile (E-modile) proposed in (4.1).
	We generate the sample size $n=10^5$ data of $Y$ from the following two distributions: Normal(0,1) and Laplace(1,2). The simulation results of sample E-modiles and absolute error (AE=$|\text{modile}-\text{E-modile}|$) based on $\tau=0.1,\ldots,0.9$ are shown in Table 2, which are based on 500 simulation replications. We can see from both Table 2 that
	the E-modile is very chose to modile in Table 1, thus the estimation method is valid.

\textcolor{blue}{
	\begin{table}[htp]
		\footnotesize
		\caption{The means and standard deviations (in parentheses) of E-modile under different errors and $\tau$s for simulation example. }
		\centering
		\begin{tabular}{@{}r|cc|cc@{}}
			\hline
			$\bm{Y}$&Normal(0,1)&&Laplace(1,2)\\
			\hline
			$\tau$&E-modile&AE&E-modile&AE\\
			\hline
			0.1& -1.097~(0.034)& 0.027~(0.020)& -1.197~(0.078)& 0.064~(0.045)\\
			0.2& -0.694~(0.026)& 0.021~(0.016)& -0.384~(0.07)& 0.055~(0.044)\\
			0.3& -0.423~(0.024)& 0.019~(0.015)& 0.156~(0.068)& 0.055~(0.040)\\
			0.4& -0.202~(0.023)& 0.019~(0.014)& 0.591~(0.068)& 0.055~(0.040)\\
			0.5& -0.001~(0.022)& 0.018~(0.013)& 0.994~(0.070)& 0.056~(0.041)\\
			0.6& 0.201~(0.022)& 0.017~(0.013)& 1.400~(0.068)& 0.055~(0.040)\\
			0.7& 0.424~(0.024)& 0.019~(0.014)& 1.844~(0.072)& 0.057~(0.043)\\
			0.8& 0.693~(0.026)& 0.021~(0.015)& 2.383~(0.069)& 0.057~(0.039)\\
			0.9& 1.097~(0.033)& 0.026~(0.019)& 3.197~(0.078)& 0.063~(0.047)\\
			\hline
		\end{tabular}
	\end{table}
}
	
	\subsection{Real data example}
	To illustrate the practical usefulness of application of our proposed methods, a daily data of S\&P500 index between January 4, 2010 and March 15, 2023 with 3321 observations in total.
	The data is downloaded from the website of Yahoo Finance (https://hk.finance.yahoo.com).
	The daily returns are computed as 100 times the difference of the log of the prices, that is, $Y_t=100\ln (p_t/p_{t-1})$, where $p_t$ is the daily price.
	Table 3 collects the summary statistics of $\{Y_t\}$, where the
	sample skewness 0.724 indicates possible asymmetries in the volatility, and the sample
	kurtosis 13.070 implies heavy tail of $\{Y_t\}$.
	Figure 4 also gives the time series plot for S\&P500.
	\par
	The sample modiles, quantiles and expectiles based on $\tau=0.1,\ldots,0.9$ are shown in Table 4 and Figure 5. The results show that
	modiles are the largest under small $\tau$ and smallest under large $\tau$
	than quantiles and expectiles. This means that modile is more conservative than
	quantiles and expectiles.
	
	\begin{table}[htp]
		\footnotesize
		\caption{Summary statistics for S\&P500 returns.  }
		\centering
		\begin{tabular}{@{}ccccccc@{}}
			\hline
			Mean&Median&std.Dev.&Skewness&Kurtosis&Min&Max\\
			\hline
			-0.037&-0.063&1.124&0.724&13.070&-8.97&12.765\\
			\hline
		\end{tabular}
	\end{table}

	\begin{table}[htp]
		\footnotesize
		\caption{The sample modile, quantile and expectile
			under different errors and $\tau$s for real data example. }
		\centering
		\begin{tabular}{@{}r|ccc@{}}
			\hline
			$\tau$&modile&quantile&expectile\\
			\hline
			0.1&-0.762&-1.177&-0.899\\
			0.2&-0.562&-0.720&-0.561\\
			0.3&-0.150&-0.439&-0.352\\
			0.4&-0.063&-0.226&-0.187\\
			0.5&0.105&0.063&-0.037\\
			0.6&0.151&0.087&0.114\\
			0.7&0.333&0.273&0.289\\
			0.8&0.358&0.576&0.523\\
			0.9&0.791&1.138&0.925\\
			\hline
		\end{tabular}
	\end{table}
	
	\begin{figure}
		\centering
		\includegraphics[width=4in]{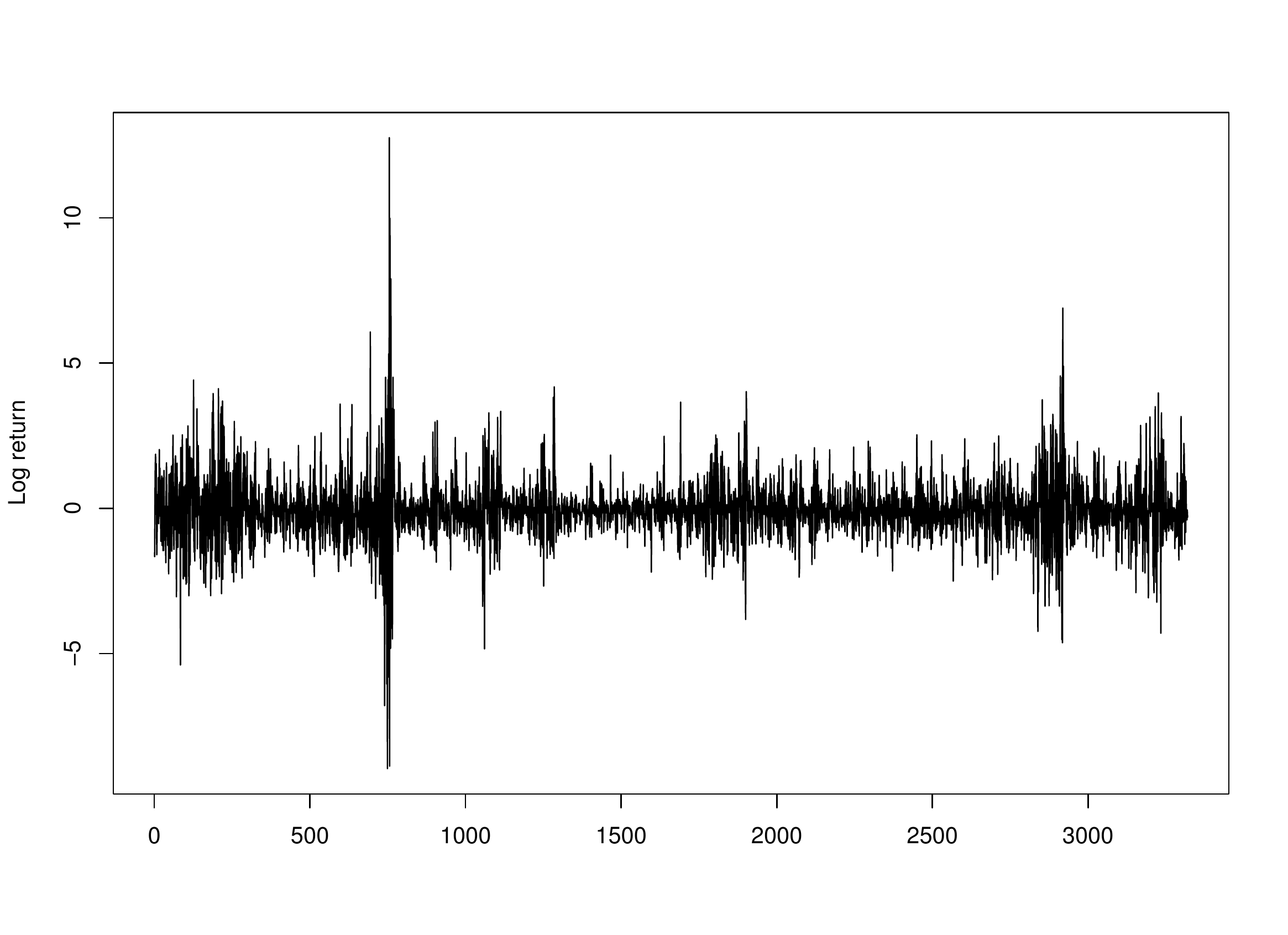}
		\caption{The time series plot of the daily return series between January 4, 2010 and March 15, 2023 .}
	\end{figure}
	\begin{figure}
		\centering
		{\label{figa} 
			\includegraphics[width=4.5in]{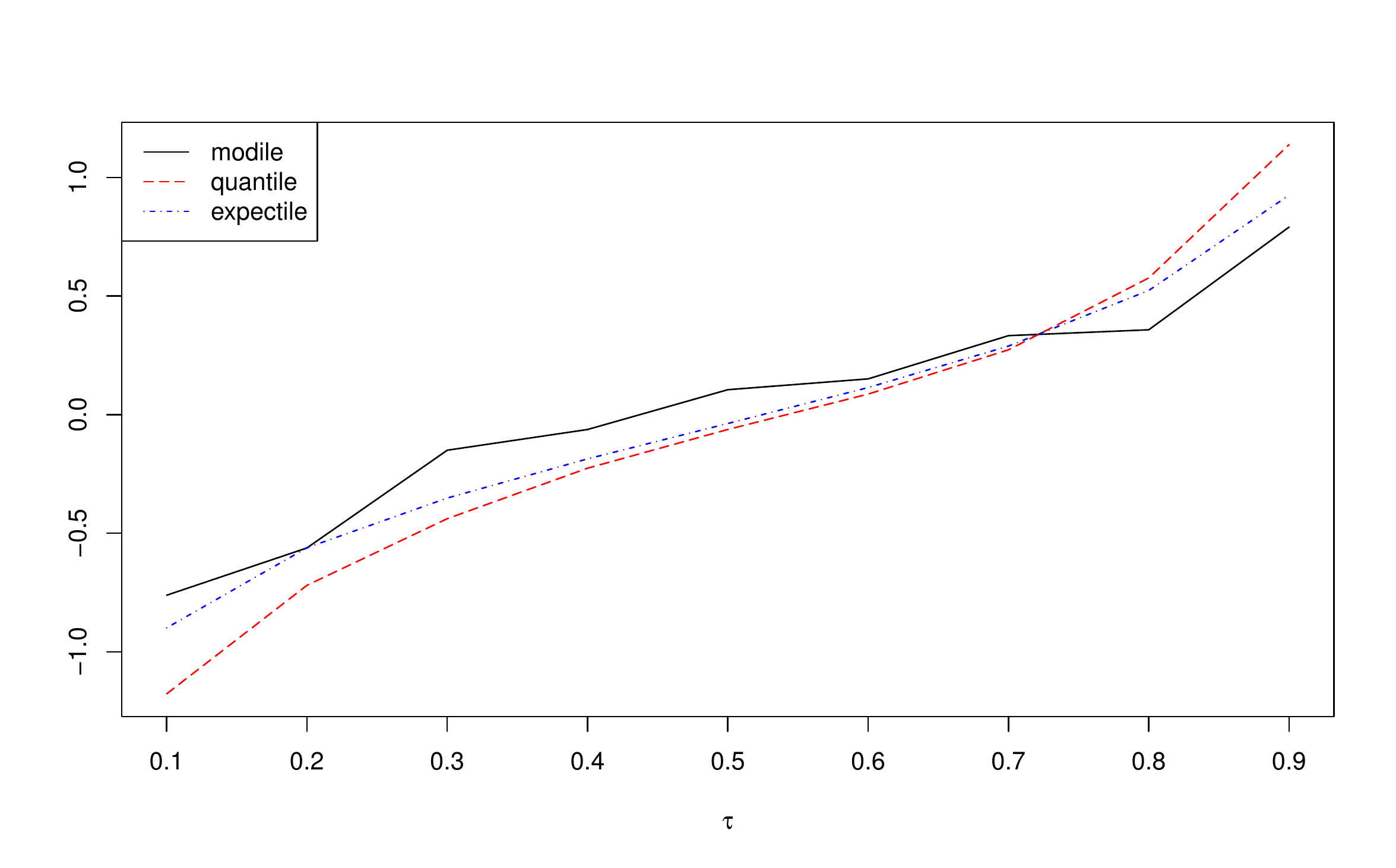}}
		\caption{The modile, quantile and expectile for real data example.}
	\end{figure}

	\section{Disussion}
	
	{\it Modiles}, as a new tail rik measure, have been introduced in this paper, which have their merits, giving existing {it quantiles} and {\it expectiles} tail risk measures.
	Like {\it quantiles}, {\it modiles} have a good interpretation and also include the unique interpretation of {\it expectiles}.
	Unlike {\it quantiles},  multivariate {\it modiles} should be  well-defined along the same line as that of univariate case. They are worth further investigation in more details.
	The applications and their deep link to current double-truncated risk measures deserve to be studied further.
	\\
	\par
	\appendix
	\section{Proof of main results} 
	{\bf Proof of Theorem 1.}
	From the condition
	$
	{f_Y(\nu_{\tau}+x)}/{f_Y(\nu_{\tau}-x)}=(1-\tau)/{\tau}$ for any $x\in(0,\delta)$, we have
	\begin{equation*}
		\begin{split}
			P(\nu_{\tau} \leq Y \leq \nu_{\tau}+\delta)
			=&\int_{\nu_{\tau}}^{\nu_{\tau}+\delta}f_Y(t)\,dt
			=\int_{0}^{\delta}f_Y(\nu_{\tau}+t)\,dt\\
			=& \frac{1-\tau}{\tau}\int_{0}^{\delta}f_Y(\nu_{\tau}-t)\,dt
			=\frac{1-\tau}{\tau}\int_{\nu_{\tau}-\delta}^{\nu_{\tau}}f_Y(t)\,dt\\
			=&\frac{1-\tau}{\tau}P(\nu_{\tau}-\delta \leq Y \leq \nu_{\tau}),
		\end{split}
	\end{equation*}
	and
	\begin{equation*}
		\begin{split}
			&E\{(Y-\nu_{\tau}) \,I(\nu_{\tau} \leq Y \leq \nu_{\tau}+\delta)\}\\
			=&\int_{\nu_{\tau}}^{\nu_{\tau}+\delta}(t-\nu_{\tau})f_Y(t)\,dt
			=\int_{0}^{\delta}tf_Y(\nu_{\tau}+t)\,dt\\
			=&\frac{1-\tau}{\tau}\int_{0}^{\delta}tf_Y(\nu_{\tau}-t)\,dt
			=\frac{1-\tau}{\tau}\int_{\nu_{\tau}-\delta}^{\nu_{\tau}}(\nu_{\tau}-t)f_Y(t)\,dt\\
			=&\frac{1-\tau}{\tau}E\{(\nu_{\tau}-Y)\,I(\nu_{\tau}-\delta \leq Y \leq \nu_{\tau})\}.
		\end{split}
	\end{equation*}
	\\
	\par
	{\bf Proof of Theorem 2.} Note that (4.1), $\hat{\nu}_{\tau}$ can be rewritten as
	\begin{equation}
		\begin{split}
			\hat{\nu}_{\tau}
			=&\arg\max_{\theta\in \mathbb{R}}\left\{\tau\hat{F}_n(\theta+h_2)-(1-\tau)\hat{F}_n(\theta-h_1)  \right\}\\
			=&\arg\max_{\theta\in \mathbb{R}}Z_n(\theta),
		\end{split}
	\end{equation}
	where $\hat{F}_n(t)=n^{-1}\sum_{i=1}^nI(Y_i\leq t)$ and
	$
	Z_n(\theta)=\tau\{\hat{F}_n(\theta+h_2)-\hat{F}_n(\nu_{\tau}+h_2)\}
	-(1-\tau)\{\hat{F}_n(\theta-h_1)-\hat{F}_n(\nu_{\tau}-h_1)\}.
	$
	We decompose $Z_n(\theta)$ into
	\begin{equation}
		\begin{split}
			Z_n(\theta)=\mu+n^{-1/2}H_n,
		\end{split}
	\end{equation}
	where $\mu=\tau\left\{F_Y(\theta+h_2)-F_Y(\nu_{\tau}+h_2)\right\}
	-(1-\tau)\left\{F_Y(\theta-h_1)-F_Y(\nu_{\tau}-h_1) \right\}$ and
	\begin{equation*}
		\begin{split}
			n^{-1/2}H_n=&\tau\left[\left\{\hat{F}_n(\theta+h_2)-F_Y(\theta+h_2)    \right\} -\left\{\hat{F}_n(\nu_{\tau}+h_2)-F_Y(\nu_{\tau}+h_2)\right\}\right]
			\\
			&-(1-\tau)\left[\left\{\hat{F}_n(\theta-h_1)-F_Y(\theta-h_1)    \right\} -\left\{\hat{F}_n(\nu_{\tau}-h_1)-F_Y(\nu_{\tau}-h_1)\right\}\right].
		\end{split}
	\end{equation*}
	Here $\mu$ is the actual deviation, while $n^{-1/2}H_n$ represents the random deviation. We can
	approximate $\mu$ by a second order Taylor expansion:
	\begin{equation}
		\begin{split}
			\mu=&\left\{\tau f_Y(\nu_{\tau}+h_2)-(1-\tau)f_Y(\nu_{\tau}-h_1)\right\}(\theta-\nu_{\tau})\\
			&+\frac{1}{2}\left\{\tau f'_Y(\nu_{\tau}+h_2)-(1-\tau)
			f'_Y(\nu_{\tau}-h_1)\right\}(\theta-\nu_{\tau})^2\{1+o(1)\}\\
			=&\frac{1}{2}\left\{\tau f'_Y(\nu_{\tau}+h_2)-(1-\tau)
			f'_Y(\nu_{\tau}-h_1)\right\}(\theta-\nu_{\tau})^2\{1+o(1)\},
		\end{split}
	\end{equation}
	where we use that $\tau f_Y(\nu_{\tau}+h_2)-(1-\tau)f_Y(\nu_{\tau}-h_1)=0$ by equation (2.3).
	\par
	Now we consider $n^{-1/2}H_n$. $n^{-1/2}H_n$ can be rewritten in two ways, depending on whether $\theta>\nu_{\tau}$ or $\theta\leq\nu_{\tau}$. If $\theta>\nu_{\tau}$,
	\begin{equation*}
		\begin{split}
			n^{-1/2}H_n=&\frac{\tau}{n}\sum_{i=1}^nI(\nu_{\tau}+h_2\leq Y_i\leq \theta+h_2)
			-\frac{1-\tau}{n}\sum_{i=1}^nI(\nu_{\tau}-h_1\leq Y_i\leq \theta-h_1)-\mu,
		\end{split}
	\end{equation*}
	and if $\theta\leq\nu_{\tau}$,
	\begin{equation*}
		\begin{split}
			n^{-1/2}H_n=&-\frac{\tau}{n}\sum_{i=1}^nI(\theta+h_2\leq Y_i\leq \nu_{\tau}+h_2)
			+\frac{1-\tau}{n}\sum_{i=1}^nI(\theta-h_1\leq Y_i\leq \nu_{\tau}-h_1)-\mu.
		\end{split}
	\end{equation*}
	We assume $\theta>\nu_{\tau}$ and $\theta\approx\nu_{\tau}$ form now on. Situation $\theta\leq\nu_{\tau}$ has a similar result. The expected value of $n^{-1/2}H_n$ is 0 and its variance is
	\begin{equation}
		\begin{split}
			Var(n^{-1/2}H_n)=&Var\left\{\frac{\tau}{n}\sum_{i=1}^nI(\nu_{\tau}+h_2\leq Y_i\leq \theta+h_2)
			-\frac{1-\tau}{n}\sum_{i=1}^nI(\nu_{\tau}-h_1\leq Y_i\leq \theta-h_1)\right\}\\
			=&\frac{\tau^2}{n}P(\nu_{\tau}+h_2\leq Y\leq \theta+h_2)
			+\frac{(1-\tau)^2}{n}P(\nu_{\tau}-h_1\leq Y\leq \theta-h_1)\\
			=&\left\{\frac{\tau^2}{n}f_Y(\nu_{\tau}+h_2)+\frac{(1-\tau)^2}{n}f_Y(\nu_{\tau}-h_1)\right\}(\theta-\nu_{\tau})
			\{1+o(1)\}\\
			=&\frac{\tau}{n}f_Y(\nu_{\tau}+h_2)(\theta-\nu_{\tau})\{1+o(1)\},
		\end{split}
	\end{equation}
	where the last equation is according to (2.3). From (A.1)-(A.4), set $t=\theta-\nu_{\tau}$, we find that
	\begin{equation}
		\begin{split}
			\hat{\nu}_{\tau}-\nu_{\tau}=\arg\max_{t\in \mathbb{R}}\left[ n^{-1/2}\sqrt{\tau f_Y(\nu_{\tau}+h_2)}W(t)-\frac{1}{2}\left\{\tau Af'_Y(\nu_{\tau}+h_2)-(1-\tau)
			f'_Y(\nu_{\tau}-h_1)\right\}t^2\right],
		\end{split}
	\end{equation}
	where $W(t)$ is a standard two-sided Brownian motion starting in 0. For any $c>0$, $c^{-1}W(c^2t)=W(t)$. Therefore, from (A.5), we have
	$$
	aW(ct)-c^2t^2=a(a/b)^{1/3}\{ W(t/c)-(t/c)^2\},
	$$
	where $a=n^{-1/2}\sqrt{\tau f_Y(\nu_{\tau}+\delta)}$, $b=1/2\left\{\tau f'_Y(\nu_{\tau}+\delta)-(1-\tau)
	f'_Y(\nu_{\tau}-\delta)\right\}$, $\delta=h_1=h_2$ and $c=(a/b)^{2/3}$.
	Therefore, $(\hat{\nu}_{\tau,\delta}-\nu_{\tau})/c$ is the limiting distribution of $\hat{z}$.
	Then, we can proof the Theorem 2.
	
	\section{R code for the numerical minimisation solution of $\nu_{\tau}$}
	\begin{verbatim}
		tau = 0.5
		y = rnorm(100)
		a=sd(y)
		b=mean(y)
		c=mean(((y-mean(y))/sd(y))^3)
		h1=a+abs(b-c)
		h2=a+abs(b+c)
		
		z = rbind(cbind(y+h1,1-tau),cbind(y-h2,-tau))
		zs = z[order(z[,1]),]
		zi = zs[which.min(cumsum(zs[,2]))+(0:1),1]
		# The objective function is constant in this interval and is minimized mean(zi)
		# Take the mid-point of the interval as the solution
	\end{verbatim}

	\footnotesize
	\bibliographystyle{elsarticle-harv}
	\bibliography{ref}

\end{document}